\documentclass[onefignum,onetabnum]{siamart250211}

\usepackage{graphicx,url}
\usepackage{amsmath,amssymb,amsopn}
\usepackage{color}
\usepackage{calc}
\usepackage{longtable}
\usepackage{multirow}
\usepackage{mathrsfs}
\usepackage{array}
\usepackage{booktabs}
\usepackage[outline]{contour}

\usepackage{tikz}
\usepackage{pgfplotstable}
\usepackage{pgfplots}
\pgfplotsset{compat=1.18}
\usepgfplotslibrary{groupplots}
\usetikzlibrary{matrix}
\usetikzlibrary{automata, external, shapes.geometric, fit}
\usetikzlibrary{arrows.meta}
\def\rv{\vec{r}}

\def\rm{\mat{r}}

\def\R{\mathbb{R}}			
\def\E {\mathbb{E}}			

\def\D{\pazocal{D}}			
\def\L{\pazocal{L}}			%
\def\D{\pazocal{D}} 			
\def\F{\pazocal{F}}

\def\G{\pazocal{G}}

\def\T{\pazocal{T}}
\def\M{\pazocal{M}}

\def\G{\pazocal{G}}

\def\Nr{\pazocal{N}}

\DeclareMathAlphabet{\pazocal}{OMS}{zplm}{m}{n}

\renewcommand{\vec}[1]{\boldsymbol{#1}}
\newcommand{\mat}[1]{\boldsymbol{\mathrm{#1}}}

\DeclareMathOperator{\tr}{tr}

\DeclareMathOperator{\cov}{cov}

\newcommand{\norm}[1]{\left\lVert#1\right\rVert}

\usepackage{lipsum}
\usepackage{amsfonts}
\usepackage{graphicx}
\usepackage{epstopdf}
\usepackage{algorithmic}
\ifpdf
  \DeclareGraphicsExtensions{.eps,.pdf,.png,.jpg}
\else
  \DeclareGraphicsExtensions{.eps}
\fi

\newsiamremark{assumption}{Assumption}
\crefname{assumption}{Assumption}{Assumptions}
\newsiamremark{remark}{Remark}
\newsiamremark{hypothesis}{Hypothesis}
\crefname{hypothesis}{Hypothesis}{Hypotheses}
\newsiamthm{claim}{Claim}
\newsiamremark{fact}{Fact}
\crefname{fact}{Fact}{Facts}

\headers{A Neural Preconditioned Newton Method}{Y. Lee, S. Liu, J. Darbon, and G. E. Karniadakis}

\title{A Neural-Operator Preconditioned Newton Method  for Accelerated Nonlinear Solvers
\thanks{Submitted to the editors
\today. \funding{This work was funded by the DARPA-DIAL grant HR00112490484.}}}

\author{
Youngkyu Lee\thanks{Division of Applied Mathematics, Brown University, Providence, RI~\email{youngkyu\_lee@brown.edu}, \email{shanqing\_liu@brown.edu}, \email{jerome\_darbon@brown.edu}, \email{george\_karniadakis@brown.edu}.}
\and
Shanqing Liu\footnotemark[2]
\and
Jerome Darbon\footnotemark[2]
\and
George Em Karnidakais\footnotemark[2]
}

\begin{document}

\maketitle

\begin{abstract}
We propose a novel neural preconditioned Newton~(NP-Newton) method for solving parametric nonlinear systems of equations.
To overcome the stagnation or instability of Newton iterations caused by unbalanced nonlinearities, we introduce a fixed-point neural operator~(FPNO) that learns the direct mapping from the current iterate to the solution by emulating fixed-point iterations.
Unlike traditional line-search or trust-region algorithms, the proposed FPNO adaptively employs negative step sizes to effectively mitigate the effects of unbalanced nonlinearities.
Through numerical experiments we demonstrate the computational efficiency and robustness of the proposed NP-Newton method across multiple real-world applications, especially for very strong nonlinearities.

\end{abstract}

\begin{keywords}
scientific machine learning, nonlinear preconditioning, neural operator, Newton's method, hybridization
\end{keywords}

\begin{MSCcodes}
90C06, 65M55, 65F08, 65F10, 68T07
\end{MSCcodes}

\section{Introduction}
\label{sec:intro}
 The numerical solution of the nonlinear system of equations is required in many engineering applications, i.e., for a given parametric nonlinear function $\F \colon \R^{n} \to \R^{n}$, we find a vector $u \in \R^{n}$ such that
\begin{equation}
    \F(u) = 0.
\label{eqn:nonlinear}
\end{equation}
Here $\F = (\F_{1}, \ldots, \F_{n})^{T}$, $\F_{i}=F_{i}(u_{1}, \ldots, u_{n})$, and $u=(u_{1}, \ldots, u_{n})^{T}$.
The Newton's method is commonly used to solve the nonlinear system of equations~\cref{eqn:nonlinear}.
Suppose $u^{(k)}$ is the current approximate solution: a new approximate solution $u^{(k+1)}$ is computed by
\begin{equation*}
    u^{(k+1)} = u^{(k)} + \lambda^{(k)}p^{(k)},
\label{eqn:newton}
\end{equation*}
where the step size $\lambda^{(k)}$ is simply set to $1$ or determined by line search algorithms~\cite{dennis1996numerical}, or the update direction $p^{(k)}$ is computed by solving the following linear system
\begin{equation}
    J^{(k)}p^{(k)}=-\F(u^{(k)}).
\label{eqn:jacobian}
\end{equation}
Here, $J^{(k)}=\F^{\prime}(u^{(k)})$ is the Jacobian of $\F(u^{(k)})$.
For simplicity, we assume that $J^{(k)}$ is nonsingular.
When $J^{(k)}$ is singular, the Moore–Penrose inverse~\cite{penrose1955generalized}, the Levenberg-Marquardt method~\cite{levenberg1944method,marquardt1963algorithm}, and the arc-length continuation method~\cite{keller1977numerical} are commonly used to find the approximated update direction.

When the initial guess is sufficiently close to the solution, the Newton's method achieves fast, second-order convergence (under appropriate assumptions).
However, obtaining such a good initial guess is generally difficult, especially when nonlinear equations have \textit{unbalanced nonlinearities}~\cite{lanzkron1996analysis}.
Here, unbalanced nonlinearities indicate that the effective step size is dominated by a subset of the total degrees of freedom, which may lead to stagnation in the nonlinear residual curve~\cite{cai1998parallel,cai2002nonlinearly,dolean2016nonlinear}.

Recently, nonlinearly preconditioned Newton methods~\cite{cai2002nonlinearly,cai2011nepin,dolean2016nonlinear,klawonn2017nonlinear,marcinkowski2005parallel} have attracted attention to address this challenge.
The main idea of nonlinear preconditioning is to replace the original nonlinear system $\F(u)=0$ with an equivalent nonlinear system $\tilde{\F}(u)=0$.
Note that two nonlinear systems are equivalent if they have the same solution.
Specifically, $\tilde{\F}$ may take the form of a composite function
\begin{equation*}
    \tilde{\F}(u):=\M(\F(u)) \text{ or } \F(\M(u)).
\end{equation*}
Here, $\M \colon \R^{n} \to \R^{n}$ is the preconditioner.
ASPIN~\cite{cai2002nonlinearly,marcinkowski2005parallel} and RASPEN~\cite{dolean2016nonlinear} utilize the additive Schwarz method or restricted additive Schwarz method as the left-preconditioner $\M$, which gives the nonlinearly preconditioned residual
\begin{equation*}
    \F_{L}(u) = \M(\F(u)) = x - \M(x).
\end{equation*}
That is, the left-preconditioning strategy first computes the residual using a fixed-point method and applies the Newton's method to solve $\F_{L}(u)=0$.
On the other hand, NEPIN~\cite{cai2011nepin}, nonlinear FETI-DP and BDDC methods~\cite{klawonn2017nonlinear} utilize the nonlinear elimination technique as the right preconditioner $\M$, which gives the nonlinearly preconditioned residual
\begin{equation*}
    \F_{R}(u) = \F(\M(u)).
\end{equation*}
Hence, the right-preconditioning strategy first applies the right-preconditioner $\M$ using a fixed-point iteration, and then applies the Newton's method to solve $\F_{R}(u)=0$.
Furthermore, several approaches that design the efficient coarse space for ASPIN and RASPEN~\cite{heinlein2022adaptive,heinlein2020additive} and combine with the field-split technique~\cite{liu2015field} have been proposed.

Recently, operator learning, a subfield of scientific machine learning, offers a powerful alternative to traditional methods for solving nonlinear system of equations.
In operator learning, a neural operator~(NO) is employed to approximate the solution operator of a given problem by learning the mapping between function spaces, with applications that include resolving a solution from material properties, determining an initial condition from a future state~(i.e., solving an inverse problem), and mapping between different state variables~\cite{li2021fourier,lu2021learning,lu2022comprehensive,wang2021learning}.
A significant advantage of this framework is that the trained NO can infer the solution for new input function in real-time without costly retraining.
However, a notable limitation is the spectral bias of neural networks, which prioritizes learning low-frequency modes.
This bias makes it difficult to capture high-frequency modes, making it challenging to obtain machine-accurate solution~\cite{rahaman2019spectral}.

On the other hand, the idea of enhancing the convergence of iterative solvers through hybridization with machine learning has recently gained significant attention.
Various approaches have been explored in this direction, such as accelerating Krylov subspace methods~\cite{pmlr-v202-kaneda23a,luo2024neural}, learning or correcting the Newton directions~\cite{NEURIPS2024_dae8afc6,jin2025fast}, developing more effective and adaptive preconditioners~\cite{gotz2018machine,heinlein2021combining,kopanivcakova2025leveraging,pmlr-v202-li23e}, and leveraging the spectral bias of neural networks to design more efficient coarse-level problems in multilevel solvers~\cite{cui2022fourier,kopanivcakova2025deeponet,lee2025fast,lee2024automatic,lee2025automatic,zhang2024blending}.
These studies demonstrate that the synergy between numerical algorithms and machine learning is gradually expanding, aiming to achieve faster convergence rates, enhanced robustness, and improved generalization capabilities across diverse problem settings.

In this work, we propose a novel nonlinearly right-preconditioning strategy to enhance the performance of the Newton's method using the pretrained neural operator.
The key idea is to mitigate the stagnation arising from unbalanced nonlinearities in a nonlinear system of equations by leveraging the neural operator.
To this end, we introduce a fixed point neural operator~(FPNO) that learns the direct mapping from the current iterate to the solution by emulating the fixed-point iteration.
By learning the underlying fixed-point map, the FPNO effectively captures the nonlinear relationship between successive Newton iterates, leading to faster and more robust convergence.
Enabling the FPNO to use negative step sizes allows the neural preconditioned Newton method to overcome stagnation or instability via negative Newton directions.
Furthermore, FPNO can employ many kinds of neural operator architectures such as DeepONets~\cite{jin2022mionet,lu2021learning,lu2022comprehensive}, FNOs~\cite{li2023fourier,li2021fourier}, Transformers~\cite{luotransolver,ovadia2024vito,shih2025transformers,wu2024Transolver}.
Numerical results confirm that the proposed neural preconditioned Newton method achieves fast convergence and reduced computational time.

The rest of the paper is organized as follows.
We first introduce the nonlinearly right-preconditioned Newton method and propose the fixed point neural operator that is used to build the neural preconditioner in~\Cref{sec:nonlinear_newton}.
\Cref{sec:benchmark} introduces the benchmark problems to demonstrate the capabilities of the proposed neural preconditioned Newton method.
The numerical experiments appear in~\Cref{sec:experiments}.
Finally, we conclude this paper with a summary in~\Cref{sec:conclusions}.
\section{A neural preconditioned Newton method}
\label{sec:nonlinear_newton}
In this section, we first briefly describe the nonlinearly right-preconditioning strategy, which is employed to construct the neural preconditioned Newton method in this paper.
Let us recall that we find a vector $u^{\ast} \in \R^{n}$ such that $\F(u^{\ast})=0$.
A nonlinear preconditioner $\M$ of the nonlinear function $\F$ is defined as $\M \approx \F^{-1}$.

The right-preconditioned Newton method applies the Newton's method to the nonlinearly right-preconditioned residual
\begin{equation*}
    \tilde{\F}(u)=\F(\M(u)).
\label{eqn:preconditioned_residual}
\end{equation*}
Starting from the initial guess $u^{(0)}$, the iteration of the right-preconditioned Newton method is defined as
\begin{equation*}
\begin{aligned}
    v^{(k)} &= \M(u^{(k)}),\\
    u^{(k+1)} &= v^{(k)} - \lambda^{(k)}(\F^{\prime}(v^{(k)}))^{-1}\F(v^{(k)}).
\end{aligned}
\end{equation*}
Here, the step size $\lambda^{(k)}$ is typically determined using a line search~\cite{dennis1996numerical}, trust-region algorithm~\cite{sorensen1982newton} or simply set to $1$.

On the other hand, in general, it is not easy to obtain such a preconditioner $\M$ explicitly, but such an operator can be defined implicitly as a fixed-point nonlinear iteration, which is given by $u = \M(u)$.
Inspired by the fixed-point nonlinear iteration, we construct the nonlinear right-preconditioner $\M$ using a neural operator $\G$ and make the neural operator learn the fixed-point nonlinear iteration.
Motivated by the fixed-point nonlinear iteration, we design the nonlinear right-preconditioner $\M$ with a neural operator $\G$, training $\G$ to emulate the fixed-point nonlinear iteration.

Let $N$ denote the input size of the neural operator.
Suppose we have a canonical restriction operator $R \colon \R^{n} \to \R^{N}$ and a canonical prolongation operator $P \colon \R^{N} \to \R^{n}$.
Then, the neural operator $\G \colon \R^{N} \to \R^{N}$ is trained to approximate the inverse of the nonlinear function $\F$ such that
\begin{equation}
    \G (v) \approx (R \circ \F^{-1} \circ P)(v), \quad \forall v \in \R^{N}.
\label{eqn:neural_preconditioner}
\end{equation}
Then, the right-preconditioner $\M \colon \R^{n} \to \R^{n}$ is naturally defined as
\begin{equation*}
    \M = P \circ \G \circ R.
\end{equation*}
Without loss of generality, we assume that $n = N$, which enables us to set $\M = \G$.

In order to emulate the fixed-point nonlinear iteration, we propose a novel neural operator called \textit{Fixed Point Neural Operator}~(FPNO) given by
\begin{equation}
\begin{aligned}
    r &= \F(u), \quad \tilde{r}=r / \norm{r}_{2},\\
    \eta &= \tanh{(\norm{r}_{2} \Nr(\tilde{r}))}, \\
    \G(u)&:=u+\eta \cdot \G_{B}(u).
\end{aligned}
\end{equation}
Here, $\Nr$ denotes a trainable neural network that takes the normalized residual vector $\tilde{r}$ and produces a scalar output, while $\G_{B}$ represents a backbone neural operator that takes the vector $u \in \R^{N}$ as the input and produces a correction vector.
When the iterate is converging to the solution, the norm of the residual $\norm{r}_{2}$ is going to $0$, which implies the step size $\eta \to 0$ and $\G(u) = u$.
That is, the proposed FPNO learns an approximate inverse of the nonlinear function $\F$ by emulating the fixed-point iteration, enabling it to act as a right-preconditioner for the Newton's method.

\begin{figure}
    \centering
    \includegraphics[width=\linewidth, trim={0cm, 8cm, 0cm, 3cm}]{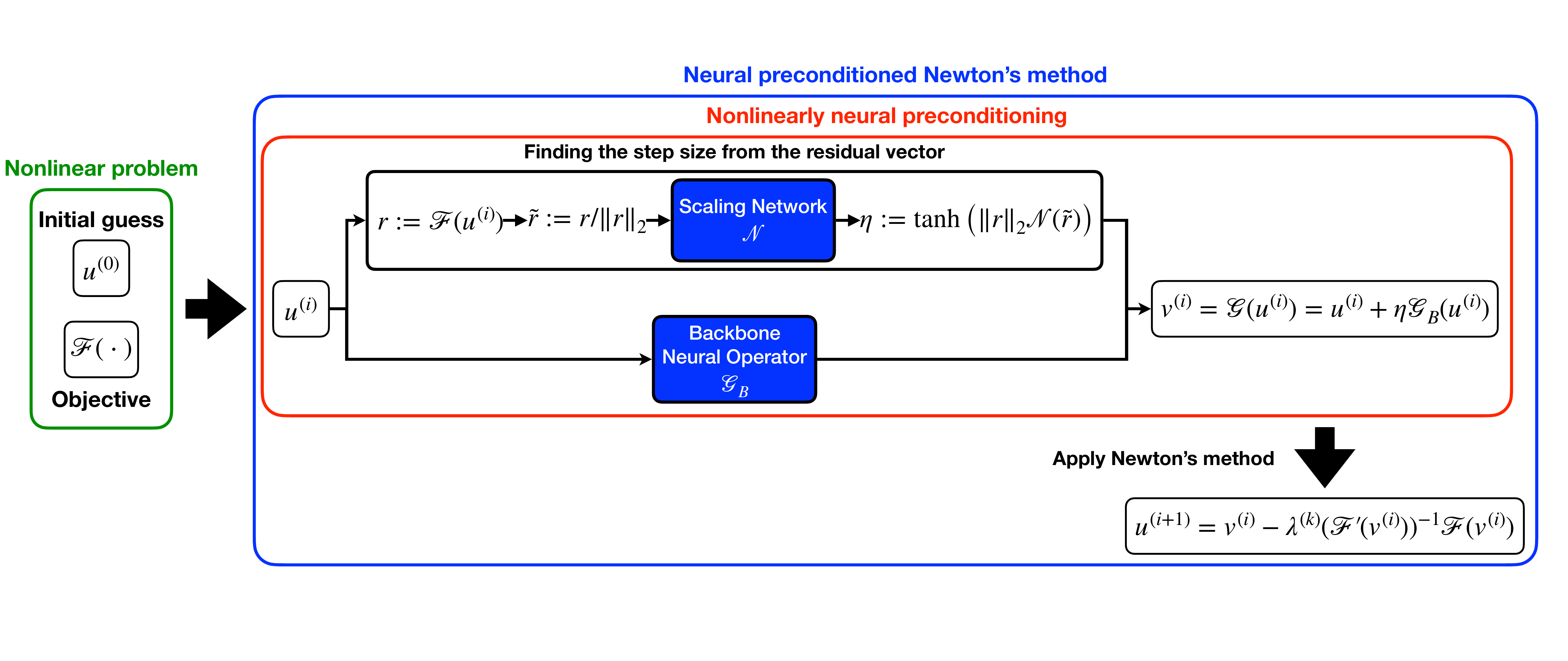}
    \caption{A schematic of neural preconditioned Newton method using the fixed point neural operator~(FPNO) with backbone neural operator $\G_{B}$.
    Starting from the initial guess $u^{(0)}$, the FPNO is first applied to get the smoothed iterate $v^{(i)}$ and the Newton's method is applied to obtain the next iterate $u^{(i+1)}$.}
    \label{fig:fpno}
\end{figure}

Furthermore, since the range of the hyperbolic tangent function is $(-1, 1)$, the predicted step size $\eta$ can be negative, which allows the FPNO to learn more robust iterations.
For example, if the Newton direction computed by~\cref{eqn:jacobian} fails to reduce the residual, the line search or the trust-region algorithm may choose an extremely small step size, such as $10^{-14}$, which can lead to stagnation.
In contrast, the FPNO can employ a negative step size and thus avoids this problem.
A graphical description of the proposed neural preconditioned Newton method using FPNO is shown in~\Cref{fig:fpno}.

\begin{remark}
The backbone neural operator $\G_{B}$ is flexible and can be implemented using various neural operator architectures, such as DeepONet~\cite{lu2021learning}, FNO~\cite{li2021fourier}, or DeepOKAN~\cite{abueidda2025deepokan}.
When solving parametric PDEs, the structure of $\G_{B}$ can be further adapted to account for the PDE parameters.
Specifically, the PDE parameter can be included as an additional input to $\G_{B}$, allowing the operator to learn parameter-dependent solution mappings.
In such cases, multi-input neural operators like MIONet~\cite{jin2022mionet} are particularly suitable.
In this work, we use MIONet as the backbone neural operator to effectively handle parametric PDEs.
\end{remark}

Lastly, in order to efficiently train the FPNO, the data generation is important.
We first sample an initial guess $u_{j}^{(0)}$ from a Gaussian random field~(GRF) and scale it to the range $10^{-4}$–$10^{-2}$.
Next, we solve $\F(u)=0$ starting from the initial guess $u_{j}^{(0)}$ and collect the Newton iteration snapshots
\begin{equation*}
    \D_{j}=\left\{\left( u_{j}^{(i)}, u_{j}^{\ast} \right)\right\}_{i=0}^{m_{j}},
\end{equation*}
where $m_{j}$ denotes the total number of iterations of the Newton method starting from the initial guess $u_{j}^{(0)}$ and $u_{j}^{\ast}$ is the reference solution.
When solving parametric PDEs, the PDE parameter $\zeta_{j}$ is sampled from a given distribution, and the corresponding Newton iteration snapshots are collected as
\begin{equation*}
    \D_{j}=\left\{\left( \zeta_{j}, u_{j}^{(i)}, u_{j}^{\ast} \right)\right\}_{i=0}^{m_{j}}.
\end{equation*}
Note that, in the case of DeepONet variants, we also collect the coordinates $x_{j}$ corresponding to the degrees of freedom.
Finally, the training dataset $\D$ is defined as
\begin{equation*}
    \mathcal{D}:=\{ \D_{j} \}_{j=1}^{m},
\end{equation*}
where $m$ denotes the total number of random initial guesses.
Note that, the total number of training samples is 
$\sum_{j=1}^{m} m_{j}$.

Finally, during the training stage, the FPNO learns the fixed-point nonlinear iteration mapping from the input $u_{j}^{(i)}$ to the reference solution $u_{j}^{\ast}$. This implies that when the input is far from $u_{j}^{\ast}$, the FPNO effectively learns a direct update from the input to the solution, whereas when the input is close to $u_{j}^{\ast}$, it learns the Newton update direction.

\section{Benchmark problems and implementation details}
\label{sec:benchmark}
In this section, we present a collection of benchmark problems for demonstrating the performance of the proposed neural preconditioned Newton method.

\subsection{Nonlinear Poisson equation}
\label{sec:nonlinear_poisson}
\begin{figure}
    \centering
    \includegraphics[width=0.5\linewidth]{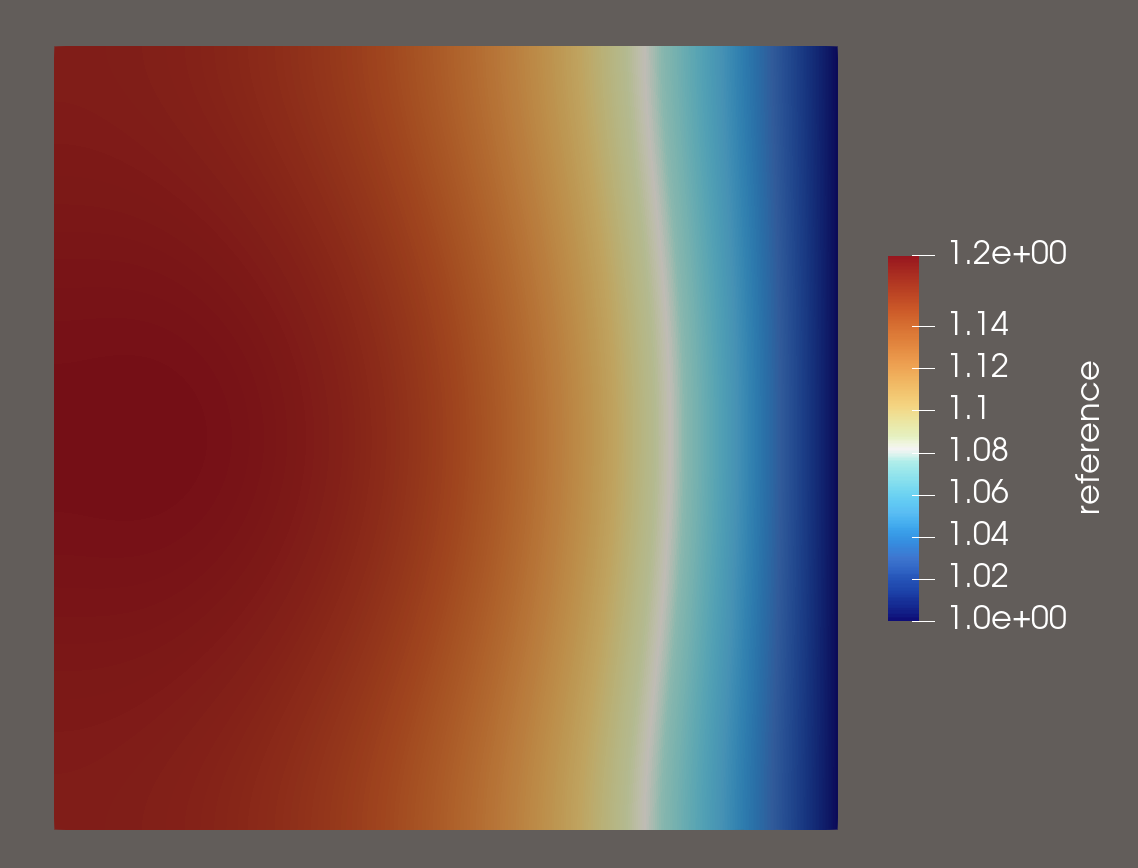}
    \caption{The numerical solution of nonlinear Poisson equation~\eqref{eqn:fem_nonlinear_poisson} computed on $129 \times 129$ mesh, where $q(u)=0.01+u^{2}$ and $f(x,y)=\sin(\pi x)\sin(\pi y)$.}
    \label{fig:nonlinear_poisson}
\end{figure}

Let $\Omega = (0, 1)^{2} \subset \mathbb{R}^{2}$ and $\Gamma =\{(x,y) \in \partial\Omega | x=1\}$.
The two-dimensional nonlinear Poisson equation is given as
\begin{equation}
\left\{\begin{aligned}
-\nabla \cdot (q(u)\nabla u) &= f, \,\,\, \text{ in } \Omega, \\
u &= 1, \text{ on } \Gamma, \\
\frac{\partial u}{\partial n} &= 0, \text{ on } \partial\Omega \setminus \Gamma,
\end{aligned}
\right.
\label{eqn:nonlinear_poisson}
\end{equation}
where $q(u)=0.01+u^{2}$ and $f$ stands for the forcing term, respectively.
We discretized the domain $\Omega$ into linear triangular elements with a mesh size of $h$, denoted as $\T_{h}$.
Then, solving the equation~\eqref{eqn:nonlinear_poisson} is equivalent to finding $u \in V$ such that
\begin{equation}
    \mathcal{F}(u)=\int_{\Omega}q(u)\nabla u \nabla v dx - \int_{\Omega} fvdx=0, \quad \forall v \in V.
    \label{eqn:fem_nonlinear_poisson}
\end{equation}
Here, $V=\{ v \in H^{1}(\Omega) \colon v\vert_{T} \in P^{1}(T) \quad \forall T \in \T_{h}, u=1 \text{ on } \Gamma \}$.
Note that $P^{1}(T)$ is the space of piecewise linear polynomials defined on $T$.
\Cref{fig:nonlinear_poisson} presents an example of the numerical solution of~\cref{eqn:fem_nonlinear_poisson} computed on $129 \times 129$ mesh.

The equation~\cref{eqn:nonlinear_poisson} is parameterized in terms of the forcing term $f$.
The forcing term $f$ is sampled using Gaussian random field~(GRF) with mean $\E[f(x)]=0.0$ and the covariance
\begin{equation}
\cov(f(x), f(y))=\sigma^{2}e^{-\frac{\norm{x-y}^{2}}{2\ell^{2}}} \quad \forall x,y \in \Omega.
\label{eqn:grf_covariance}
\end{equation}
The parameter $\sigma$ and $\ell$ are chosen as $\sigma=0.1$ and $\ell=0.1$.

\subsection{Hyper elasticity}
\label{sec:hyper_elasticity}
We next consider the compressible Neo-Hookean model, a popular hyper elasticity model, applied to a square with an oval-shaped hole.
Let $\Omega \subset (0,1)^{2} \subset \R^{2}$ and $\partial \Omega = \Gamma_{1} \cup \Gamma_{2} \cup \Gamma_{3}$, where $\Gamma_{1}$, $\Gamma_{2}$, and $\Gamma_{3}$ are the bottom boundary, top boundary, and other boundaries, respectively.
We discretized the domain $\Omega$ into linear quadrilateral elements with a mesh size of $h$, denoted as $\T_{h}$.
Then, the minimization problem of the two-dimensional compressible Neo-Hookean model is given as
\begin{equation}
\min_{u \in V} \mathcal{L}(u):= \int_{\Omega} \frac{\mu}{2}(I_{c}-3)-\mu \ln(J_{F})+\frac{\lambda}{2}(J_{F}-1)^{2}dx,
\label{eqn:hyper_elasticity}
\end{equation}
where 
\begin{equation*}
\left\{\begin{aligned}
V &= \{ u \in [H^{1}(\Omega)]^{2} \colon u\vert_{T} \in [P^{1}(T)]^{2} \quad \forall T \in \T_{h}, \\ &u=(0,0) \text{ on } \Gamma_{1} \text{ and } u=(0, u_{t}) \text{ on } \Gamma_{2} \}, \\
\mu &= \frac{E}{2(1+\nu)}, \lambda= \frac{E\nu}{(1+\nu)(1-2\nu)}, \\
F &= I+\nabla u, I_{c} = \tr(F^{T}F), J_{F}=\det(F). \\
\end{aligned}
\right.
\end{equation*}
Note that we set $E=1.0$ and $\nu=0.49$ to model rubber-like materials, which withstand large deformations.

\begin{figure}
    \centering
    \includegraphics[width=0.5\linewidth]{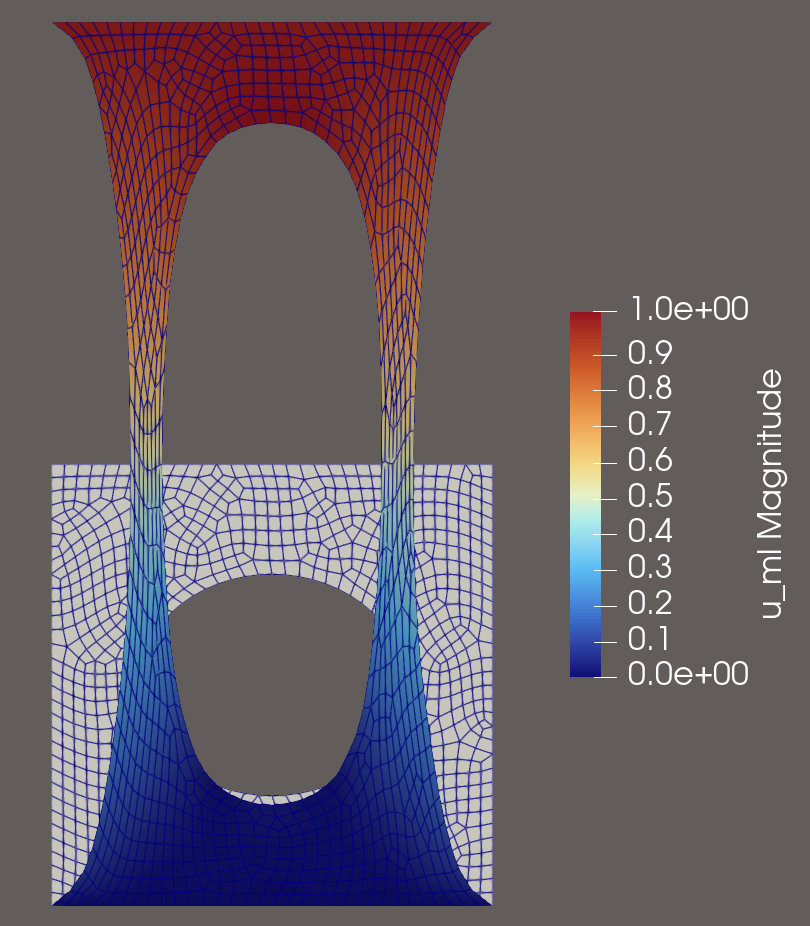}
    \caption{The computed displacement of hyper elasticity problem~\cref{eqn:hyper_elasticity} in $\Omega \subset (0,1)^{2}$ when $u_{t}=1$ and the maximum mesh size $h = 1/32$. Note that the gray object is the initial mesh.}
    \label{fig:square_hole_mesh}
\end{figure}

The minimization problem~\eqref{eqn:hyper_elasticity} is solved by finding a solution $u \in V$ that satisfies the optimality condition
\begin{equation}
    \frac{\partial \mathcal{L}}{\partial u}=\mathcal{F}(u)=0.
    \label{eqn:fem_hyper_elasticity}
\end{equation}
\Cref{fig:square_hole_mesh} shows an example of the computed displacement $u$ and the initial mesh.
The problem~\cref{eqn:hyper_elasticity} is parameterized in terms of the $y$-displacement $u_{t}$ on the top boundary, which is sampled from the uniform distribution $U(0, 2)$.

\subsection{Implementation details}
\label{sec:implementation}
We use the FEniCS library~\cite{BarattaEtal2023} to implement the finite element discretization of the benchmark problems and to generate the training samples for the FPNOs.
To handle parametric PDEs, we adopt the MIONets~\cite{jin2022mionet} as the backbone neural operator of the FPNOs and refer to the resulting model as the FP-MIONets.
The FP-MIONets are implemented using PyTorch~\cite{paszke2019hh} and trained using the AdamW optimizer~\cite{ilya2019decoupled}, with a batch size of $100$, a weight decay of $5 \times 10^{-4}$, and a learning rate of $10^{-4}$.
The training process stops if the average relative $L^{2}$ error of the validation samples does not improve for 1,000 consecutive epochs.
Note that the relative mean squared error loss function is used for training, which is given as
\begin{equation*}
	\L_{\text{rel}}(u_{\theta}, u_{\text{ref}}) = \frac{1}{m}\sum_{i=1}^{m}\frac{(u^{(i)}_{\theta}-u^{(i)}_{\text{ref}})^{2}}{(u^{(i)}_{\text{ref}})^{2}+\varepsilon},
\end{equation*}
where $u_{\theta}=(u^{(1)}_{\theta}, \ldots, u^{(m)}_{\theta})$ and $u_{\text{ref}}=(u^{(1)}_{\text{ref}}, \ldots, u^{(m)}_{\text{ref}})$ are the prediction of the neural network and the reference solution, respectively.
Note that the symbol $m$ denotes the batch size, and a small constant $\varepsilon=10^{-4}$ is added to prevent the loss term from blowing up when $u^{(i)}_{\text{ref}} = 0$.
Details about the network architectures, size of datasets, and training times are described in~\Cref{sec:further_details}.

The proposed neural preconditioned Newton methods are implemented using PETSc library~\cite{balay2019petsc}.
The preconditioning of the Newton's method with the FP-MIONet is implemented using the petsc4py interface.
Note that the python implementation is not compiled and conducted without hyperparameter tuning, suggesting that the practical performance can be further enhanced.
All numerical experiments are conducted using the Oscar supercomputer at Brown University, whose computing node is equipped with an AMD EPYC 9554 64-Core Processor~(256GB) and an NVIDIA L40S GPU~(48GB).
\section{Numerical results}
\label{sec:experiments}
In this section, we present numerical experiments to evaluate the performance of the proposed neural preconditioned Newton method.
The Newton method with cubic backtracking line search~(Newton-LS) and the Newton method with trust region~(Newton-TR) are chosen as the baseline method.
We construct the neural preconditioned Newton method combining the pretrained FP-MIONets with Newton-LS and Newton-TR, denoted as NP-Newton-LS and NP-Newton-TR, respectively.
All local problems associated with the computation of the Newton direction~\cref{eqn:jacobian} are solved by the MUMPS direct solver~\cite{amestoy2001fully}.
All Newton methods terminate when either of the following criteria is satisfied:
\begin{equation*}
	\Vert \rv^{(i)}\Vert_{2} \leq 10^{-15} \quad \text{ or } \quad \frac{\Vert \rv^{(i)}\Vert_{2}}{\Vert \rv^{(0)}\Vert_{2}} \leq 10^{-9}.
\end{equation*}

\subsection{Nonlinear Poisson equation with various forcing term}
First, we solve the nonlinear Poisson equation~\cref{eqn:nonlinear_poisson} with Newton-LS and NP-Newton-LS.
We consider the following three cases of the forcing term $f$ to verify the robustness of the proposed NP-Newton-LS:
\begin{itemize}
    \item I: $f \equiv 1$.
    \item II: $f \sim GRF(\nu, \sigma, \ell)$, where $\nu=0.0$, $\sigma=0.1$, and $\ell=0.1$. Note that the FP-MIONet is trained on this case.
    \item III: $f \sim GRF(\nu, \sigma, \ell)$, where $\nu=0.0$, $\sigma=1.0$, and $\ell=0.1$. Unlike the previous case, $f$ can take the high frequency value.
\end{itemize}
Here, $GRF(\nu, \sigma, \ell)$ denotes the GRF with mean $\nu$ and the covariance defined in~\cref{eqn:grf_covariance} using $\sigma$ and $\ell$.
Furthermore, to demonstrate robustness in resolution and highlight the effectiveness of the proposed NP-Newton-LS particularly on fine meshes, we conduct tests on two meshes $\T_{h_{1}}$ and $\T_{h_{2}}$. Note that $h_{1}=1/32$ and $h_{2}=1/128$, where the FP-MIONet is mainly trained on the coarse mesh $\T_{h_{1}}$.

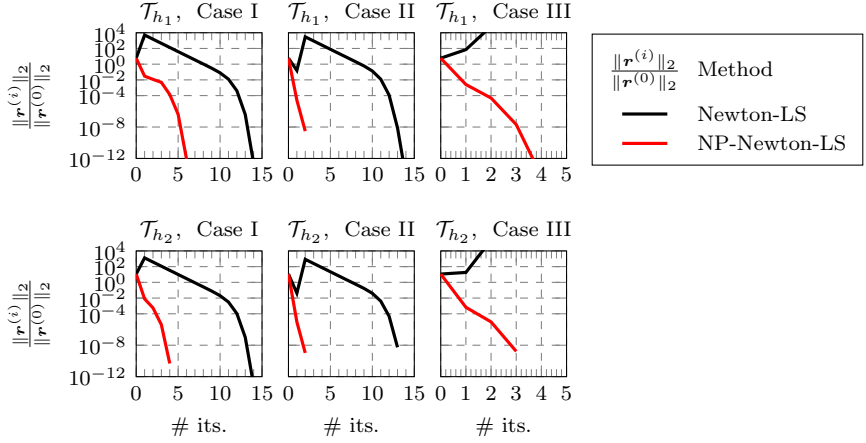
\begin{figure}
\centering
\tikzexternaldisable
\tikzsetnextfilename{nonlinear_poisson}
 \begin{tikzpicture}
   \begin{groupplot}[
       group style={
           group size = 3 by 2,
           horizontal sep = 10pt,
           vertical sep = 35pt,
         },
       legend pos=north east,
       width=0.25\textwidth,
       height=0.25\textwidth,
       minor x tick num=4,
       grid=major, 
       grid style={dashed,gray}, 
       xmode=normal,
       ymode=log,
	   ytick={1e4, 1e2, 1e0, 1e-2, 1e-4, 1e-8, 1e-12},
       xtick={0, 5, 10, 15},
       ymin=1e-12,          
       ymax=1e4,        
       xmin=0,
       xmax=15,
       title style={yshift=-1.5ex,font=\footnotesize},           
       tick label style={font=\footnotesize},
       label style={font=\footnotesize},
     ]

     \nextgroupplot[align=left,
       title={$\T_{h_{1}}$, \ Case I},
       ylabel= {$\frac{\| \boldsymbol{r}^{(i)} \|_{2}}{\|\boldsymbol{r}^{(0)}  \|_{2}}$},
       xlabel= {},            
       legend style={at={(0.6, 0.77)},anchor=west}
     ]

     \addplot[color = black, very thick] table [x index=0, y=Newton-LS, col sep=comma] {results/np_coarse/case_1.csv};
     \addplot[color = red, very thick] table [x index=0, y=Hybrid-Newton-LS, col sep=comma] {results/np_coarse/case_1.csv};

     \nextgroupplot[align=left,
       title={$\T_{h_{1}}$, \ Case II},
       yticklabels={},
       xlabel= {},                       
       ylabel= {},          
        xlabel= {},                            
       legend style={at={(0.6, 0.77)},anchor=west}
     ]
     
     \addplot[color = black, very thick] table [x index=0, y=Newton-LS, col sep=comma] {results/np_coarse/case_2.csv};
     \addplot[color = red, very thick] table [x index=0, y=Hybrid-Newton-LS, col sep=comma] {results/np_coarse/case_2.csv};

     \nextgroupplot[align=left,
       title={$\T_{h_{1}}$, \ Case III}, 
       yticklabels={},
       xlabel= {},                          
       ylabel= {},
       xmax = 5,
       xtick= {0, 1, 2, 3, 4, 5},
       legend style={at={(0.6, 0.77)},anchor=west}
     ]

     \addplot[color = black, very thick] table [x index=0, y=Newton-LS, col sep=comma] {results/np_coarse/case_3.csv};
     \addplot[color = red, very thick] table [x index=0, y=Hybrid-Newton-LS, col sep=comma] {results/np_coarse/case_3.csv};

     \nextgroupplot[align=left,
       title={$\T_{h_{2}}$, \ Case I},
       ylabel= {$\frac{\| \boldsymbol{r}^{(i)} \|_{2}}{\|\boldsymbol{r}^{(0)}  \|_{2}}$},
       xlabel= {\# its.},           
       legend style={at={(0.6, 0.77)},anchor=west}
     ]
     
     \addplot[color = black, very thick] table [x index=0, y=Newton-LS, col sep=comma] {results/np_fine/case_1_fine.csv};
     \addplot[color = red, very thick] table [x index=0, y=Hybrid-Newton-LS, col sep=comma] {results/np_fine/case_1_fine.csv};
    
     \nextgroupplot[align=left,
       title={$\T_{h_{2}}$, \ Case II},      
       yticklabels={},
       ylabel= {},          
       xlabel= {\# its.},                   
       legend style={at={(0.6, 0.77)},anchor=west}
     ]

     \addplot[color = black, very thick] table [x index=0, y=Newton-LS, col sep=comma] {results/np_fine/case_2_fine.csv};
     \addplot[color = red, very thick] table [x index=0, y=Hybrid-Newton-LS, col sep=comma] {results/np_fine/case_2_fine.csv};

     \nextgroupplot[align=left,
       yticklabels={},
       title={$\T_{h_{2}}$, \ Case III},              
       xlabel= {\# its.},       
       ylabel= {},      
       xmax = 5,
       xtick= {0, 1, 2, 3, 4, 5},
       legend style={at={(0.6, 0.77)},anchor=west}
     ]
     
     \addplot[color = black, very thick] table [x index=0, y=Newton-LS, col sep=comma] {results/np_fine/case_3_fine.csv};
     \label{pgfplot:newton_ls}
     \addplot[color = red, very thick] table [x index=0, y=Hybrid-Newton-LS, col sep=comma] {results/np_fine/case_3_fine.csv};
     \label{pgfplot:hybrid_newton_ls}
                            
\end{groupplot}
\matrix [ draw, matrix of nodes, anchor = north, node font=\footnotesize,
   column 1/.style={nodes={align=center,text width=1.0cm}},
   column 2/.style={nodes={align=left,text width=2.0cm}},        
   ] at ($(group c2r1) + (5.0, 0.825)$)
   {
$\frac{\| \boldsymbol{r}^{(i)} \|_{2}}{\|\boldsymbol{r}^{(0)}  \|_{2}}$ & Method	  \\
\ref{pgfplot:newton_ls} & Newton-LS \\ 
\ref{pgfplot:hybrid_newton_ls}  & NP-Newton-LS\\
};
\end{tikzpicture}
\caption{The convergence of the Newton's method for the nonlinear Poisson equation. The first row shows the results for the coarse mesh, and the second row shows the results for the fine mesh. In Case III, the Newton-LS method diverges, and the iteration is terminated when the relative residual norm exceeds $10^{4}$.}
\label{fig:convergence_nonlinear_poisson}
\end{figure}

\Cref{fig:convergence_nonlinear_poisson} compares the convergence behavior of the Newton-LS method and the proposed NP-Newton-LS method for the nonlinear Poisson equation.
For the NP-Newton-LS, the FP-MIONet that is pretrained in the coarse mesh $\T_{h_{1}}$ is employed as a nonlinear right-preconditioner to enhance the convergence of Newton’s method.
As shown in the figure, the Newton-LS occasionally suffers from an increase in the residual norm during the Newton iteration, particularly in Case III, where the method eventually diverges.
This phenomenon is mainly attributed to the unbalanced nonlinearities of given nonlinear problem, which lead to an unstable search direction in the standard Newton iteration.
In contrast, the NP-Newton-LS effectively mitigates this issue, maintaining a monotonic decrease in the residual norm across all cases, demonstrating the stabilizing effect of the proposed neural preconditioning strategy.

\Cref{tab:nonlinear_poisson,tab:nonlinear_poisson_fine} summarize the number of iterations and the computational times for the coarse mesh $\T_{h_{1}}$ and the fine mesh $\T_{h_{2}}$, respectively.
In all test cases, the NP-Newton-LS significantly reduces the number of iterations compared to the Newton-LS, achieving substantial computational speed-up.
Notably, in Cases I and II, the iteration counts are reduced by more than 50\%, and the total computational time decreases accordingly.
In Case III, the Newton-LS fails to converge due to the unbalanced nonlinearities, whereas the NP-Newton-LS successfully converges within only a few iterations.
These results clearly demonstrate that the neural operator based nonlinear preconditioning improves both the convergence rate and the robustness of Newton’s method.

\begin{table}
\centering
\caption{The number of iterations and the computational time in seconds~(s) required by the Newton-LS and NP-Newton-LS for different cases in the coarse mesh $\T_{h_{1}}$. In Case III, since the relative residual norm of Newton-LS exceeds $10^{4}$, the Newton iteration is terminated and noted as divergent.}
\label{tab:nonlinear_poisson}
\begin{tabular}{@{}clccc@{}}
\toprule
Case & Method & \# iterations & Time~(s) & Speed-up \\ \midrule
\multirow{2}{*}{I} & Newton-LS & 14 & 0.0710 & -  \\ \addlinespace
 & NP-Newton-LS & 6 & 0.0671 & 5.81\% \\
\midrule\midrule
\multirow{2}{*}{II} & Newton-LS & 14 & 0.0795 & -  \\ \addlinespace
 & NP-Newton-LS & 2 & 0.0336 & 136.60\%\\
 \midrule\midrule
\multirow{2}{*}{III} & Newton-LS & diverge & - & - \\ \addlinespace
 & NP-Newton-LS & 4 & 0.0455 & $\infty$ \\
 \bottomrule
\end{tabular}
\end{table}

\begin{table}
\centering
\caption{The number of iterations and the computational time in seconds~(s) required by the Newton-LS and NP-Newton-LS for different cases in the fine mesh $\T_{h_{2}}$. Similar to the coarse mesh, in Case III, since the relative residual norm of Newton-LS exceeds $10^{4}$, the Newton iteration is terminated and noted as divergent.}
\label{tab:nonlinear_poisson_fine}
\begin{tabular}{@{}clccc@{}}
\toprule
Case & Method & \# iterations & Time~(s) & Speed-up \\ \midrule
\multirow{2}{*}{I} & Newton-LS & 14 & 0.5435 & -  \\ \addlinespace
 & NP-Newton-LS & 4 & 0.2768 & 96.35\% \\
\midrule\midrule
\multirow{2}{*}{II} & Newton-LS & 13 & 0.5030 & -  \\ \addlinespace
 & NP-Newton-LS & 2 & 0.1575 & 219.37\%\\
 \midrule\midrule
\multirow{2}{*}{III} & Newton-LS & diverge & - & - \\ \addlinespace
 & NP-Newton-LS & 3 & 0.1985 & $\infty$ \\
 \bottomrule
\end{tabular}
\end{table}

\subsection{Hyper elasticity problem with small and large deformation}

Next, we solve the hyper elasticity problem~\cref{eqn:hyper_elasticity} with Newton-LS, Newton-TR, NP-Newton-LS, and NP-Newton-TR. 
In order to highlight the performance of the proposed NP-Newton-LS and NP-Newton-TR, we consider the following two cases of the top displacement corresponding to small and large deformations:

\begin{itemize}
    \item I: $u_{t}= 0.1$. In this case, the computed displacement is quite similar with the displacement computed from the linear elasticity model.
    \item II: $u_{t}=1$. In this case, the Newton–LS method fails to converge without additional techniques.
    To address this issue, we employ the incremental loading technique~\cite{ogden1984non}, a common approach for mitigating unbalanced nonlinearities, with a loading step size of $\delta u_{t}=0.1$.
    Specifically, the total loading is applied gradually in several increments rather than in a single step, which allows the Newton's method to trace the equilibrium path more stably.
    At each increment, the solution obtained from the previous step serves as an initial guess for the next, thereby improving convergence and preventing divergence due to large deformation.
    The resulting method is referred to as IC–Newton–LS.
\end{itemize}

Moreover, we also consider the coarse mesh $\T_{h_{1}}$ and the fine mesh $\T_{h_{2}}$ to further investigate the robustness of the proposed methods with respect to mesh resolution.
The corresponding numbers of degrees of freedom for $\T_{h_{1}}$ and $\T_{h_{2}}$ are summarized in~\Cref{tab:dofs}.

\begin{figure}
\centering
\tikzexternaldisable
\tikzsetnextfilename{hyper_elasticity}
 \begin{tikzpicture}
   \begin{groupplot}[
       group style={
           group size = 2 by 2,
           horizontal sep = 10pt,
           vertical sep = 35pt,
         },
       legend pos=north east,
       width=0.35\textwidth,
       height=0.25\textwidth,
       minor x tick num=4,
       grid=major, 
       grid style={dashed,gray}, 
       xmode=normal,
       ymode=log,
	   ytick={1e4, 1e2, 1e0, 1e-2, 1e-4, 1e-8, 1e-12},
       xtick={0, 5, 10, 15},
       ymin=1e-12,          
       ymax=1e4,        
       xmin=0,
       xmax=15,
       title style={yshift=-1.5ex,font=\footnotesize},           
       tick label style={font=\footnotesize},
       label style={font=\footnotesize},
     ]

     \nextgroupplot[align=left,
       title={$\T_{h_{1}}$, \ Case I},
       ylabel= {$\frac{\| \boldsymbol{r}^{(i)} \|_{2}}{\|\boldsymbol{r}^{(0)}  \|_{2}}$},
       xlabel= {},            
       legend style={at={(0.6, 0.77)},anchor=west}
     ]

     \addplot[color = black, very thick] table [x index=0, y=Newton-LS, col sep=comma] {results/hy_coarse/e_1_g_0.1.csv};
     \addplot[color = red, very thick, dashed] table [x index=0, y=Hybrid-Newton-LS, col sep=comma] {results/hy_coarse/e_1_g_0.1.csv};
     \addplot[color = blue, very thick] table [x index=0, y=Newton-TR, col sep=comma] {results/hy_coarse/e_1_g_0.1.csv};
     \addplot[color = green, thick] table [x index=0, y=Hybrid-Newton-TR, col sep=comma] {results/hy_coarse/e_1_g_0.1.csv};

     \nextgroupplot[align=left,
       title={$\T_{h_{1}}$, \ Case II},
       yticklabels={},
       xlabel= {},                       
       ylabel= {},
       xmax = 110,
       xtick = {0, 20, 40, 60, 80, 110},
       legend style={at={(0.6, 0.77)},anchor=west}
     ]
     
     \addplot[color = black, very thick, dotted] table [x index=0, y=IC-Newton-LS, col sep=comma] {results/hy_coarse/e_1_g_1.csv};
     \addplot[color = blue, very thick] table [x index=0, y=Newton-TR, col sep=comma] {results/hy_coarse/e_1_g_1.csv};
     \addplot[color = green, very thick] table [x index=0, y=Hybrid-Newton-TR, col sep=comma] {results/hy_coarse/e_1_g_1.csv};

     \nextgroupplot[align=left,
       title={$\T_{h_{2}}$, \ Case I},
       ylabel= {$\frac{\| \boldsymbol{r}^{(i)} \|_{2}}{\|\boldsymbol{r}^{(0)}  \|_{2}}$},
       xlabel= {\# its.},   
       xmax = 25,
       xtick = {0, 5, 10, 15, 20, 25},
       legend style={at={(0.6, 0.77)},anchor=west}
     ]
     
     \addplot[color = black, very thick] table [x index=0, y=Newton-LS, col sep=comma] {results/hy_fine/e_1_g_0.1_fine.csv};
     \label{pgfplot:hy_newton_ls}
     \addplot[color = red, very thick, dashed] table [x index=0, y=Hybrid-Newton-LS, col sep=comma] {results/hy_fine/e_1_g_0.1_fine.csv};
     \label{pgfplot:hy_hybrid_newton_ls}
     \addplot[color = blue, very thick] table [x index=0, y=Newton-TR, col sep=comma] {results/hy_fine/e_1_g_0.1_fine.csv};
     \addplot[color = green, thick] table [x index=0, y=Hybrid-Newton-TR, col sep=comma] {results/hy_fine/e_1_g_0.1_fine.csv};
    
     \nextgroupplot[align=left,
       title={$\T_{h_{2}}$, \ Case II},      
       yticklabels={},
       ylabel= {},          
       xlabel= {\# its.},     
       xmax = 220,
       xtick = {0, 40, 80, 120, 160, 220},
       legend style={at={(0.6, 0.77)},anchor=west}
     ]

     \addplot[color = black, very thick, dotted] table [x index=0, y=IC-Newton-LS, col sep=comma] {results/hy_fine/e_1_g_1_fine.csv};
     \label{pgfplot:hy_ic_newton_ls}
     \addplot[color = blue, very thick] table [x index=0, y=Newton-TR, col sep=comma] {results/hy_fine/e_1_g_1_fine.csv};
     \label{pgfplot:hy_newton_tr}
     \addplot[color = green, very thick] table [x index=0, y=Hybrid-Newton-TR, col sep=comma] {results/hy_fine/e_1_g_1_fine.csv};
     \label{pgfplot:hy_hybrid_newton_tr}
                            
\end{groupplot}
\matrix [ draw, matrix of nodes, anchor = north, node font=\footnotesize,
   column 1/.style={nodes={align=center,text width=1.0cm}},
   column 2/.style={nodes={align=left,text width=2.2cm}},        
   ] at ($(group c2r1) + (4.0, 0.825)$)
   {
$\frac{\| \boldsymbol{r}^{(i)} \|_{2}}{\|\boldsymbol{r}^{(0)}  \|_{2}}$ & Method	  \\
\ref{pgfplot:hy_newton_ls} & Newton-LS \\ 
\ref{pgfplot:hy_hybrid_newton_ls}  & NP-Newton-LS\\
\ref{pgfplot:hy_ic_newton_ls}  & IC-Newton-LS\\
\ref{pgfplot:hy_newton_tr} & Newton-TR \\ 
\ref{pgfplot:hy_hybrid_newton_tr}  & NP-Newton-TR\\
};
\end{tikzpicture}
\caption{Convergence of Newton's method for the hyper elasticity problem. Since the IC-Newton-LS uses the incremental loading technique, the history of the relative residual norm is concatenated. In Case I, the convergence behaviors of NP-Newton-LS and NP-Newton-TR are almost identical, so their plots overlap.}
\label{fig:convergence_hyper_elasticity}
\end{figure}
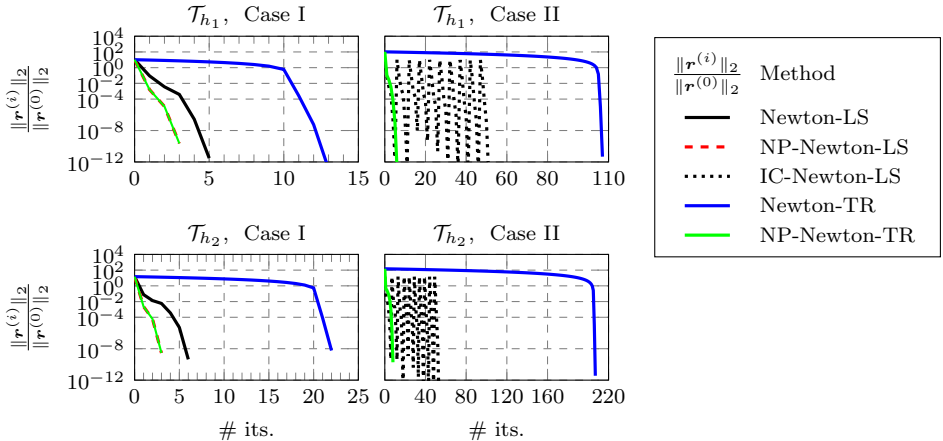

\Cref{fig:convergence_hyper_elasticity} provides the convergence histories for the relative residual norm on two different meshes $\T_{h_{1}}$ and $\T_{h_{2}}$.
The performance of the proposed neural preconditioned Newton methods~(NP-Newton-LS and NP-Newton-TR) is compared against the Newton's methods~(Newton-LS and Newton-TR) and the incremental loading method~(IC-Newton-LS).
Similar to the case of the nonlinear Poisson equation, the FP-MIONet used to construct the neural preconditioner is pretrained in the coarse mesh $\T_{h_{1}}$.
In the simpler Case I, all methods successfully converge.
The proposed NP-Newton-LS and NP-Newton-TR exhibit the fastest convergence in terms of iteration count, requiring only $3$ iterations on both meshes.
Note that the convergence behaviors of NP-Newton-LS and NP-Newton-TR are nearly identical, and their plots overlap.
The effectiveness and robustness of the proposed framework are most evident in the more challenging Case II.
Since the Newton-LS method fails to converge, the Newton–TR method is used as the baseline when the incremental loading technique is not applied.
It requires over $100$ and $200$ iterations for the $\T_{h_{1}}$ and $\T_{h_{2}}$ meshes, respectively.
By employing the incremental loading technique, the IC-Newton-LS method significantly outperforms the standard Newton-TR method in terms of stability and convergence, with convergence attained in roughly $40$ iterations.
On the other hand, the proposed NP-Newton-TR demonstrates exceptional performance.
It converges robustly in fewer than $10$ iterations on both coarse and fine meshes.

\Cref{tab:hyper_elasticity,tab:hyper_elasticity_fine} summarize the number of iterations and the computational times for both coarse and fine meshes.
In Case I, although the NP-Newton-LS reduces the number of iterations by approximately half compared with the Newton-LS, the actual computational speed-up is negative for both coarse and fine meshes. However, the degradation is noticeably smaller on the fine mesh than on the coarse mesh.
On the other hand, the NP-Newton-TR significantly reduces the number of iterations compared with the Newton-TR, achieving actual computational speed-ups of $71.03\%$ and $246.47\%$, respectively.
In Case I where the deformation is minimal, the Newton's method does not suffer from unbalanced nonlinearity.
Therefore, since the cost of preconditioning using a neural operator is higher, we can conclude that no practical speedup can be achieved.

However, the situation changes drastically under the large deformation.
In Case II, due to the unbalanced nonlinearities, the Newton-LS fails to converge entirely, while the Newton-TR requires a substantially large number of iterations to reach a solution.
A well-known approach to mitigate this issue is the incremental loading technique, denoted as IC-Newton-LS.
As reported in both tables, IC-Newton-LS effectively stabilizes the solution process, significantly reducing the total number of iterations and achieving actual computational speed-ups of $143.04\%$ and $331.33\%$ on the coarse and fine meshes, respectively, compared with the Newton-TR.
In contrast, our proposed NP-Newton-TR exhibits remarkable efficiency and robustness: it converges to the solution in only $6$ iterations on the coarse mesh and $8$ iterations on the fine mesh.
This corresponds to extraordinary actual computational speed-ups of $597.89\%$ and $1130.46\%$, respectively, demonstrating the superiority of the NP-Newton-TR in handling hyper elasticity problem with large deformation.
These results highlight that the combination of the proposed FPNO-based neural preconditioning and the trust-region algorithm provides a highly robust and efficient solver for challenging nonlinear problems.

\begin{table}
\centering
\caption{The number of iterations and the computational time in seconds~(s) required by the Newton-LS, Newton-TR, IC-Newton-LS, NP-Newton-LS, and NP-Newton-LS for different cases in the coarse mesh $\T_{h_{1}}$.}
\label{tab:hyper_elasticity}
\begin{tabular}{@{}clccc@{}}
\toprule
Case & Method & \# iterations & Time~(s) & Speed-up \\ \midrule
\multirow{4}{*}{I} & Newton-LS & 5 & 0.0596 & -  \\ 
\addlinespace
 & NP-Newton-LS & 3 & 0.0901 & -33.85\% \\
 \addlinespace
& Newton-TR & 13 & 0.1529 & - \\
 \addlinespace
 & NP-Newton-TR & 3 & 0.0894 & 71.03\% \\
\midrule\midrule
\multirow{3}{*}{II} & Newton-TR & 107 & 1.2541 & - \\
\addlinespace
 & IC-Newton-LS & 42 & 0.5160 & 143.04\% \\
 \addlinespace
 & NP-Newton-TR & 6 & 0.1797 & 597.89\% \\
 \bottomrule
\end{tabular}
\end{table}

\begin{table}
\centering
\caption{The number of iterations and the computational time in seconds~(s) required by the Newton-LS, Newton-TR, IC-Newton-LS, NP-Newton-LS, and NP-Newton-LS for different cases in the fine mesh $\T_{h_{2}}$.}
\label{tab:hyper_elasticity_fine}
\begin{tabular}{@{}clccc@{}}
\toprule
Case & Method & \# iterations & Time~(s) & Speed-up \\ \midrule
\multirow{4}{*}{I} & Newton-LS & 6 & 0.2114 & -  \\ 
\addlinespace
 & NP-Newton-LS & 3 & 0.2258 & -6.38\% \\
 \addlinespace
& Newton-TR & 22 & 0.7553 & - \\
 \addlinespace
 & NP-Newton-TR & 3 & 0.2180 & 246.47\% \\
\midrule\midrule
\multirow{3}{*}{II} & Newton-TR & 207 & 6.9841 & - \\
\addlinespace
 & IC-Newton-LS & 44 & 1.6192 & 331.33\% \\
 \addlinespace
 & NP-Newton-TR & 8 & 0.5676 & 1130.46\% \\
 \bottomrule
\end{tabular}
\end{table}

\section{Summary}
\label{sec:conclusions}

We proposed a novel nonlinearly right-preconditioning framework to accelerate Newton’s method using a pretrained neural operator.
In particular, we introduced the fixed-point neural operator~(FPNO), which learns the direct mapping from the current Newton iterate to the reference solution by emulating the fixed-point iteration of a given nonlinear function.
In order to address unbalanced nonlinearities, the proposed FPNO adaptively allows the negative step sizes based on the residual direction.
Moreover, FPNO is compatible with various neural operator architectures such as DeepONet, FNO, and Transolver; in this work, MIONet was employed for implementation to handle the parametric nonlinear PDEs.
The performance of the proposed neural preconditioned Newton method was evaluated through benchmark problems, showing its capability to enhance convergence.

The presented framework can be extended in several promising directions.
To handle large-scale domains, the domain decomposition method based nonlinear preconditioners such as ASPIN~\cite{cai2002nonlinearly}, RASPEN~\cite{dolean2016nonlinear}, NEPIN~\cite{cai2011nepin}, and nonlinear FETI-DP and BDDC methods~\cite{klawonn2017nonlinear} are widely used as state-of-the-art preconditioners for the Newton's method.
These methods often rely on a coarse-level correction to further accelerate convergence; however, constructing such a structure, e.g., via the full approximation scheme
, is often nontrivial~\cite{heinlein2022adaptive,heinlein2020additive}.
This challenge could be addressed by employing the proposed neural preconditioning framework to construct efficient two-level nonlinear solvers, potentially combining data-driven learning with multilevel strategies.

\appendix
\section{Details of training procedure of FP-MIONets} 
\label{sec:further_details}
In this section, we provide details on the network architectures of the neural operators and their training strategies for all problems discussed in~\Cref{sec:benchmark}.
We employ the squeeze-and-excitation~(SE) operation~\cite{hu8578843} with a softmax activation function instead of the sigmoid activation function, and integrate it with a residual neural network~(ResNet)~\cite{he7780459} composed of fully connected layers, denoted as SE-ResNet.
Specifically, each layer of the SE-ResNet consists of a fully connected layer followed by a skip connection and an SE operation.
Note that we employ the Gaussian error linear unit~(GELU) as the activation function in all networks.
To predict vector-valued functions, the branch network of the FP-MIONet generates additional coefficients equal to the dimension of the vector, and each component of the vector-valued function is obtained by taking the inner product of these coefficients with the output of the trunk network.
The detailed architectures of the FP-MIONets are summarized in~\Cref{tab:architectures}.

Regarding the training datasets, we first summarize the discretized meshes used in the finite element method in~\Cref{tab:dofs}.
For each problem, the mesh $\T_{h_{1}}$ is employed to generate the training datasets for the FP-MIONets.
\Cref{tab:datasets} provides the training times required for all FP-MIONets, along with the number of samples used in each dataset.
Note that the relative $L^{2}$ errors on the validation sets are $0.6\%$ and $0.1\%$ for the nonlinear Poisson equation and the hyper elasticity problem, respectively, demonstrating the high accuracy of the trained FP-MIONets.

\begin{table}
\centering
\caption{The summary of FP-MIONets' architectures. Here, NP and HE denote the nonlinear Poisson and the hyper elasticity, respectively. The symbol $S$, $B$, $B_{f}$ and $T$ represent the scaling, branch, feature-branch and trunk networks in FP-MIONets, respectively. The feature-branch network is introduced to handle physical parameters and corresponds to the second branch network in the MIONet architecture~\cite{jin2022mionet}.
For the NP, the source term $f$ is used as input, whereas for the HE, the $y$-component of the displacement is utilized.
A fully connected ResNet is employed with the GELU activation function, and a squeeze-and-excitation~(SE) block is incorporated to further enhance network performance.
The symbol $[\cdot]$ indicates the width for the neural network.}
\label{tab:architectures}
\begin{tabular}{@{}lclc@{}}
\toprule
 Problem & Subnetworks & \multicolumn{1}{c}{Layers} & Activation            \\ \midrule
 \multirow{4}{*}{NP} & $S$  & SE-ResNet[1089, 512, 512, 1] & \multirow{4}{*}{GELU} \\ \addlinespace
 & $B$  & SE-ResNet[1089, 512, 512, 512, 256] &  \\ \addlinespace
 & $B_{f}$  & SE-ResNet[1089, 512, 512, 512, 256] &  \\ \addlinespace
 & $T$  & ResNet[2, 512, 512, 512, 256] &  \\ \addlinespace
\midrule\midrule
\multirow{4}{*}{HE} & $S$  & SE-ResNet[1029, 512, 512, 1] & \multirow{4}{*}{GELU} \\ \addlinespace
 & $B$  & SE-ResNet[1029, 512, 512, 512, 512] &  \\ \addlinespace
 & $B_{f}$  & SE-ResNet[1, 512, 512, 512, 256] &  \\ \addlinespace
 & $T$  & ResNet[2, 512, 512, 512, 256] &  \\
\bottomrule
\end{tabular}
\end{table}

\begin{table}
\centering
\caption{Summary of the number of degrees of freedom associated with meshes for problems in~\Cref{sec:benchmark}. Here, NP and HE denote the nonlinear Poisson and the hyper elasticity, respectively.
For the NP problem, the unit square domain is discretized using triangular elements, whereas for the HE problem, an unstructured mesh of quadrilateral elements is used to discretize a unit square containing an elliptical hole at the center with a major axis of $0.6$ and a minor axis of $0.5$.
}
\label{tab:dofs}
\begin{tabular}{@{}lcc@{}}
\toprule
Problem & $\T_{h_{1}}$ & $\T_{h_{2}}$ \\ \midrule
NP & 1,089 & 16,641 \\ \midrule
HE & 1,029 & 3,932 \\
\bottomrule
\end{tabular}
\end{table}

\begin{table}
\centering
\caption{Summary of the number of initial guesses~$(N_{g})$, generated training/validation samples~$(N_{s_{t}}/N_{s_{v}})$, the number of epochs, and training time in hours~(h) for all problems. Here, NP and HE denote the nonlinear Poisson and the hyper elasticity, respectively.}
\label{tab:datasets}
\begin{tabular}{@{}lcccc@{}}
\toprule
Problem & $N_{g}$ & $N_{s_{t}}$/$N_{s_{v}}$ & Epochs & Time~(h) \\ \midrule
NP & 3,150 & 10,658/1,185 & 1,921 & 0.62 \\ \midrule
HE & 2,000 & 46,257/5,140 & 3,293 & 4.12 \\
\bottomrule
\end{tabular}
\end{table}

\bibliographystyle{siamplain}
\bibliography{references}
\end{document}